\documentclass[reqno]{amsart}
\theoremstyle{plain}
\newtheorem{thm}{Theorem}
 
 \newtheorem{prop}[thm]{Proposition}

\newtheorem{Main}{Main Theorem}

\numberwithin{equation}{section}

\begin{document}
\large
%Topmatter
\title[Eigenvalue monotonicity]
{Eigenvalue monotonicity for the Ricci-Hamilton flow}
\author{ Li MA}

\address{Li Ma: Department of Mathematical Sciences \\
Tsinghua University \\
Beijing 100084 \\
China}

\email{lma@math.tsinghua.edu.cn}

 \dedicatory{In Memory of S.S.Chern}
\date{}
\thanks{}


\begin{abstract}In this short note, we discuss the monotonicity of the
eigen-values of the Laplacian operator to the Ricci-Hamilton flow on
a compact or a  complete non-compact Riemannian manifold. We show
that the eigenvalue of the Lapacian operator on a compact domain
associated with the evolving Ricci flow is non-decreasing provided
the scalar curvature having a non-negative lower bound and Einstein
tensor being not too negative. This result will be useful in the
study of blow-up models of the Ricci-Hamilton flow.

\end{abstract}
\maketitle
%end topmatter

We study the monotonicity property of the eigenvalues of the
Laplacian operator $\Delta:=\Delta_{g(t)}$ of the evolving metric
$(g(t))$  along the Ricci-Hamilton flow. Our main result is stated
in the Main Theorem below.

We recall the definition of Ricci-Hamilton flow. Let $(M,g_0)$ be a
compact or a complete non-compact Riemannian manifold. We are given
a family of Riemannian metrics
  $\{g(t)\}$ with $g(0)=g_0, 0\leq t<T$. Let $g=g(t)$.
In the local coordinates $(x^i)$, we write
$$
g=g_{ij}dx^idx^j.
$$
Then the Ricci-Hamilton flow is defined by
  the evolution equation for Riemannian metrics:
  $$
\partial_tg_{ij}=-2R_{ij},\;\;\;\mbox{on $M_T:=M\times [0,T)$}
  $$
where $R_{ij}$ is the Ricci tensor of the metric $g:=g(t)$ and $T$
is the maximal existing time for the flow. In \cite{H1}, R.Hamilton
proved the local existence of the flow for the compact manifold
case. His argument is much more simplified by De Turck \cite{DeT}.
When $(M, g_0)$ is a complete non-compact Riemannian manifold with
bounded geometry, W.X.Shi \cite{Sh} obtained the local existence
result for the flow. In \cite{P02}, G.Perelman introduced two
important entropy functionals which are monotone along the
Ricci-Hamitlon flow coupled with back-ward heat flows.

We consider the change rate of first eigenvalue of Laplacian
operator associated with the Ricci flow $g(t)$ in a compact domain
in $M$. More precisely, let $D$ be a compact domain with smooth
boundary in the manifold $M$.
 Let $f:=f(x,t)$ be the first eigen-function of $\Delta$. Then we
 have
 $$
-\Delta f=\mu f,\;\;\;\mbox{in $D$,}\eqno{(1)}
 $$
 with Dirichlet boundary condition
 $$
f=0,\;\;\;\mbox{on ${\partial D}$}
 $$
where $\mu:=\lambda_1(D)>0$ is the corresponding first eigenvalue of
$\Delta$. We can normalize $f$ such that
$$
\int_D f^2 dv=1.\eqno{(2)}
$$
Here $dv$ is the volume form of the metric $g=g(t)$.

We will use the following conventions: For $F=F(t)$ being a quantity
on $D$, we let $F'$ be the derivative with respect to the time
variable $t$ and let $F_i, F_{ij}, F_{ij;k}$, etc, be the covariant
derivatives of $F$ with respect to the Levi-Civita connection of the
metric $g$. We also assume that when $D=M$ is compact, we impose the
natural condition that $\int_M fdv=0$. We abuse the upper or lower
indexes in this paper when their meanings are clear. We write by
$D_T=D\times [0,T)$.

Differentiating $(1)$, we get
$$
-\Delta' f-\Delta f'=\mu' f+\mu f'.\;\;\;\mbox{on $D$,}\eqno{(3)}
$$
Multiplying both sides of $(3)$ by $f$ and integrating over $D$,
we have
$$
-\int \Delta' f\cdot f-\int \Delta f'\cdot f=\mu'+\mu\int f'f.
\eqno{(4)}
$$
Note that
$$
-f'\Delta f=\mu f'f,\;\;\;\;\mbox{on $D$.}
$$
Then we have
$$
\mu'=-\int \Delta' f f
$$

Hence using (7) in the appendix we obtain that
$$
\mu'=-2\int R_{ij}f_{ij}f.
$$
Hereafter, we abuse the use of upper and lower indices by using
orthonormal moving frame. If $n=2$, then we have
$$
R_{ij}=\frac{1}{2}Rg_{ij},
$$
and $$ \mu'=-\int R\Delta f\cdot f=\mu\int f^2Rdv\geq 0,
$$
provided $R>0$.

In the case when $R\geq C>0$ in $D\times [t_0,t)$ for some uniform
constant $C$, we have
$$
(\log\mu)'\geq C,
$$
and $\mu(t)\geq \mu(t_0)e^{C(t-t_0)}$ for $t>t_0$.

It is also clear that when $R\leq0$, we have
$$
\mu'=-\int R\Delta f\cdot f=\mu\int f^2Rdv\leq 0.
$$
If we further have that $R\leq -C\leq 0$ in $D\times [t_0,t)$ for
some uniform constant $C>0$, we have
$$
(\log\mu)'\leq -C,
$$
and $\mu(t)\leq \mu(t_0)e^{-C(t-t_0)}$ for $t>t_0$. So we have the
following result
\begin{prop} Assume $n=2$.

1. If $R\geq C>0$ in $D\times [t_0,t)$ for some uniform constant $C$
along the Ricci flow, we have
$$
(\log\mu)'\geq C,
$$
and $\mu(t)\geq \mu(t_0)e^{C(t-t_0)}$ for $t>t_0$.

2. If we have that $R\leq -C\leq 0$ in $D\times [t_0,t)$ for some
uniform constant $C>0$ along the Ricci flow, we have
$$
(\log\mu)'\leq -C,
$$
and $\mu(t)\leq \mu(t_0)e^{-C(t-t_0)}$ for $t>t_0$.
\end{prop}

We now consider the higher dimensional case when $n\geq 3$. We
denote by $E_{ij}$ the Einstein tensor
$$
E_{ij}:=R_{ij}-\frac{R}{2}g_{ij}
$$
of the metric $g=\{g_{ij}\}$.

 When $n\geq 3$, we notice that
\begin{eqnarray*}
\mu'/2=-\int R_{ij} f_{ij}f&=&\int (R_{ij}f)_jf_i\\
&=&\int R_{ij;j}ff_i-\int R_{ij}f_if_j\\
&=&\frac{1}{2}\int R_if_if+\int R_{ij}f_if_j\\
&=&-\frac{1}{2}\int R(f_if)_i+\int R_{ij}f_if_j\\
&=&-\frac{1}{2}\int R\Delta f f-\frac{1}{2}\int R|\nabla f|^2+\int R_{ij}f_if_j\\
&=&\frac{\mu}{2}\int Rf^2+\int (R_{ij}-\frac{R}{2}g_{ij})f_if_j\\
&=&\frac{\mu}{2}\int Rf^2+\int E_{ij}f_if_j.
\end{eqnarray*}
Here we have used again the contracted second Bianchi identity
$$
2R_{ij;j}=R_i.
$$

Assume that ${R}-2a\geq 0$ on $D\times\{t\}$ and
$$ E_{ij}\geq -ag_{ij},\;\;\mbox{in $D\times\{t\}$}
$$
for some constant $a$. Then we have
\begin{align*}
\mu'/2=&\frac{\mu}{2}\int Rf^2+\int E_{ij}f_if_j \\
&\geq \frac{\mu}{2}\int Rf^2-a\int |\nabla f|^2\\
&=\mu\int (\frac{R}{2}-a)f^2\geq 0.
\end{align*}

From the proof, one can see that the same result is also true for
higher order eigenvalues.
 Therefore, we conclude the following result:

\begin{Main} Let $g=g(t)$ be the evolving metric along the Ricci-Hamilton flow
with $g(0)=g_0$ being the initial metric in $M$. Let $D$ be a smooth
bounded domain in $(M, g_0)$. Let $\mu>0$ be the first eigenvalue of
the Laplacian operator of the metric $g(t)$. If there is a constant
$a$ such that the scalar curvature $R\geq 2a$ in $D\times \{t\}$ and
the Einstein tensor
$$
E_{ij}\geq -ag_{ij},\;\;\mbox{in $D\times\{t\}$,}
$$
then we have $\mu'\geq 0$, that is $\mu$ is non-increasing in $t$,
furthermore, $\mu'(t)>0$ for the scalar curvature $R$ not being the
constant $2a$. The same result is also true for higher order
eigenvalues.
\end{Main}

We remark that these result may be useful in the study of blow-up
models of Ricci-Hamilton flow on a complete Riemannian manifold
$(M,g_0)$.

\section*{Appendix}

Given a Riemannian manifold $(M^n,g)$. In the local coordinates
$(x^i)$, we write
$$
g=g_{ij}dx^idx^j.
$$

Let $u$ be a smooth function on $M$. Then the Laplacian of $u$ is
defined by
$$
\Delta u=\frac{1}{\sqrt{|g|}}\partial_i(\sqrt{|g|}g^{ij}\partial_j,
u)
$$
where
$$
(g^{ij})=(g_{ij})^{-1}
$$
is the inverse of the matrix $(g_{ij}) $ and $|g|=det(g_{ij})$.

For convenient of readers, we review some facts about
Ricci-Hamilton flow. Recall that the metric $g=g(t)$ satisfies the
Ricci-Hamilton flow:
$$
\partial_tg_{ij}=-2R_{ij},\;\;\;\mbox{on $M_T.$}
  $$
Along this flow, we have that $$
\partial_tg^{ij}=2R^{ij},\;\;\;\mbox{on $M_T.$}
$$
and
$$
\partial_tdv=-Rdv,\;\;\;\;\eqno{(5)}
$$
where $dv$ is the volume element and $R$ is the scalar curvature of
the metric $g(t)$ respectively. For smooth functions $u$ and $v$
with compact support in $M$, we have by a use of the divergence
theorem that
$$
\int g^{ij}u_iv_jdv=-\int \Delta u vdv.
$$
We write by $u^i=g^{ij}u_j$. Differentiating both sides of this
equation and using $(5)$ we have
$$
2\int R^{ij}u_iv_j-\int Ru^iv_i=-\int \Delta'uv+\int \Delta u\cdot
vR.\;\;\;\;\eqno{(6)}
$$
Since
$$
-\int Ru^iv_i=\int (Ru^i)_iv=\int R_iu^iv+\int Rv\Delta u,
$$
we have that
$$
2\int R^{ij}u_iv_j+\int \left(\nabla R,\nabla u\right)v=-\int
\Delta'uv.
$$
As before, using integration by part and using the contracted
second Bianchi identity, we obtain that
$$
2\int R^{ij}u_iv_j=-\int\left(\nabla R,\nabla u\right)v-2\int
R^{ij}u_{ij}v.
$$
Combining this with (6) we get
$$
\int \Delta'uv=2\int R^{ij}u_{ij}v.
$$
Hence we have the variation formula for the Laplacian operator to
the Ricci-Hamilton flow:
$$
\Delta'u=2R^{ij}u_{ij}.\eqno{(7)}
$$

In particular when $n=2$, we have that
$$
R_{ij}=\frac{1}{2}Rg_{ij},
$$
and
$$
\Delta'u=R\Delta u.
$$

For more material about Ricci-Hamilton flow, one may see \cite{P02}
and \cite{M}. An earlier version of this paper has some misprint in
sign, and it is in arxiv.math.DG/0403065. Using our idea, D.Kokotov
and D.Korotkin [Normalized Ricci Flow on Riemann Surfaces and
Determinant of Laplacian, Letters in Mathematical Physics,
71(3)(2005)241-242] can give a simple proof of the fact that the
determinant of Laplace operator in a smooth metric over compact
Riemann surfaces of an arbitrary genus g monotonously grows under
the normalized Ricci flow. Together with results of Hamilton and
B.Chow that under the action of the normalized Ricci flow a smooth
metric tends asymptotically to the metric of constant curvature,
this leads to a simple proof of the Osgood-Phillips-Sarnak theorem
stating that within the class of smooth metrics with fixed conformal
class and fixed volume the determinant of the Laplace operator is
maximal on the metric of constant curvature.

\section*{Acknowledgements} The author would like to thank
Prof.J. Eichhorn and other faculty there very much for the
encouragements and helpful discussions when he was visiting
Greifswald Univerisity, Germany, in the summer of 2004.

\end{document}